\documentclass[reqno]{amsart}

\usepackage{amssymb}

\newtheorem{theorem}{Theorem}

\newtheorem{corollary}{Corollary}[theorem]

\newtheorem{conjecture}{Conjecture}

\newcommand{\bnum}   {\ensuremath{\mathcal B}-number}
\newcommand{\bprime} {\ensuremath{\mathcal B}-prime}
\newcommand{\brep}   {\ensuremath{\mathcal B}-representation}

\newcommand{\congr}[3]{{\ensuremath{{#1} \equiv {#2}
      \pmod{#3}}}}
\newcommand{\bq}[2]{\ensuremath{{#1}^2 + {#1}{#2} + {#2}^2}}

\newcommand{\bqz}[1]{{\ensuremath{{#1}^2 + {#1} + 1}}}

\newcommand{\bqm}[2]{\ensuremath{{#1}^2 - {#1}{#2} + {#2}^2}}

\newcommand{\bqab}{\bq{a}{b}}
\newcommand{\bqcd}{\bq{c}{d}}
\newcommand{\bqxy}{\bq{x}{y}}
\newcommand{\bqalbet}{\bq{\alpha}{\beta}}

\newcommand{\bqmab}{\bqm{a}{b}}

\newcommand{\bqmxy}{\bqm{x}{y}}

\newcommand{\Integer}{\ensuremath{\mathbb{Z}}}
\newcommand{\Positive}{\ensuremath{\mathbb{Z^+}}}
\newcommand{\Nonnegative}{\ensuremath{\mathbb{Z^*}}}
\newcommand{\Rational}{\ensuremath{\mathbb{Q}}}
\newcommand{\NonnegativeRational}{\ensuremath{\mathbb{Q^*}}}
\newcommand{\Real}{\ensuremath{\mathbb{R}}}

\begin{document}
\title{Elementary results on the binary quadratic form  \boldmath$a^2+ab+b^2$\unboldmath}
\author{Umesh P. Nair}
\address{Mentor Graphics Corporation, 8005
  SW Boeckman Road, Wilsonville, OR 97070,
  USA.} 
\email{umesh\_nair@mentor.com}

\begin{abstract}
  This paper examines with elementary proofs some
  interesting properties of numbers in the binary
  quadratic form $a^2+ab+b^2$, where $a$ and $b$ are
  non-negative integers.  Key findings of this paper
  are (i) a prime number $p$ can be represented as
  $a^2+ab+b^2$ if and only if $p$ is of the form $6k+1$,
  with the only exception of $3$, (ii) any positive integer
  can be represented as $a^2+ab+b^2$ if and only if its
  all prime factors that are not in the same form have even
  exponents in the standard factorization, and (iii)
  all the factors of an integer in the form $a^2+ab+b^2$,
  where $a$ and $b$ are positive and relatively prime
  to each other, are also of the same form.  A
  general formula for the number of distinct
  representations of any positive integer in this form
  is conjectured.  A comparison of the results with the
  properties of some other binary quadratic forms
  is given.
\end{abstract}
  \maketitle

\section{Background}

For more than three centuries, binary quadratic forms
and their prime representations have been studied quite
extensively.  Lagrange was the first to give a complete
treatment of the topic, and various mathematicians,
including Legendre, Euler and Gauss, contributed to the
theory~\cite{DAV1999,DIC1923}.

In addition to the general theory of binary quadratic
forms, attempts were made to study the properties of
individual forms.  The form $a^2 \pm kb^2$ got
particular attention, and the prime representations of
many such forms have been determined.  A list of such
representations can be found in ~\cite[p. 71]{BER1994} and
~\cite{MWPrimeRep}.  A detailed account of primes of these
forms is given in \cite{COX1989}.

This paper investigates numbers of the form \bqab,
where $a$ and $b$ are non-negative integers, using
elementary number theory. 

\section{Notations and Definitions}
\label{sec:Notations}
The following symbols are used in this paper.

\vspace{10pt}
\begin{tabular}{lcl}
  \Integer &:& Set of integers :  \{\ldots, -2, -1, 0,
  1, 2, \ldots\} \\
  \Positive &:& Set of positive integers : \{1, 2,
  \ldots\} \\
  \Nonnegative &:& Set of non-negative integers : 
  \{0, 1, 2, \ldots\} \\
  \Rational &:& Set of rational numbers \\
  \NonnegativeRational &:& Set of non-negative rational numbers \\
  \Real &:& Set of real numbers
\end{tabular}
\vspace{10pt}

For convenience, the following definitions are used:
\begin{enumerate}
\item A {\em \brep} is the form \bqab{}, where $a,b \in
  \Nonnegative$ and $a \ge b$.  
\item Two \brep{}s \bqab{} and \bqcd{} are {\em
    distinct} if either $a \ne c$ or $b \ne d$.
\item An integer is a {\em \bnum{}} if it has at least one \brep.
\item If a \bnum{} is a prime, it is a {\em \bprime}.
\item A positive integer is said to be {\em square-free} if
  it does not have a square factor greater than $1$.  In
  other words, all of its prime factors occur only once in the
  factorization.
\end{enumerate}

\section{General results}
\subsection{Some trivial results}
The following are quite obvious from the definition, or can
be easily deduced.
\begin{theorem}
  \label{th:BqNotNegative}
  A number in the form $a^2 \pm ab + b^2$, with $a,b \in \Real$,
  is never negative.
\end{theorem}
\begin{proof}
  WLOG, assume $|a| \ge |b|$.  This means $a^2 \ge |ab|$.
  Since both $a^2 - |ab|$ and $b^2$ are always non-negative,
  so is their sum.   This proves that both \bqab{} and
  $a^2-ab+b^2$ are non-negative.
\end{proof}

\begin{theorem}
  \label{th:UniqueB}
  Given $p=\bqab$, where $a,b \in \Nonnegative$, and the
  values of $p$ and $a$, there is a unique $b$.
\end{theorem}
\begin{proof}
  Given $a$ and $p$, the value of $b$ is given by 
  \begin{equation*}
    b = \frac{-a \pm \sqrt{4p-3a^2}}{2}
  \end{equation*}
Since $b$ is nonnegative, only the value corresponding to
the positive sign applies here.
\end{proof}
\begin{theorem}
  \label{th:BNumOfAllIntegers}
  If an integer $n$ is in the form \bqab{}, where
  $a,b \in \Integer$, then $n$ is a \bnum{}, i.e.,
  $n=\bqcd$ for some $c, d \in \Nonnegative$.
\end{theorem}
\begin{proof}
  There are three cases to consider:
  \begin{enumerate}
  \item {\em If both $a$ and $b$ are non-negative:} In this case,
    $n$ is a \bnum{} by definition.
  \item {\em If both $a$ and $b$ are negative:} Set $c = -a, d =
    -b$, and now $n = \bqcd$, where $c,d \in \Nonnegative$.
  \item {\em If one of $a$ and $b$ is negative and other
      non-negative:} WLOG, assume that $a \ge 0$ and $b < 0$.  If
    $a > -b$, set $c = a+b, d=-b$; else set $c = -(a+b),
    d=a$. Now, $n = \bqcd$, where $c,d \in \Nonnegative$.
  \end{enumerate}
  Hence the result.
\end{proof}

It is to be noted that, by Theorem~\ref{th:BNumOfAllIntegers}, $a$ and $b$ in the
\brep{} can be any number in \Integer{} in order to get a
\bnum{}, but we are restricting the definition to
\Nonnegative{} only.

\subsection{Identities}

The following identities can be
verified easily.
\begin{subequations}
  \begin{align}
    c^2(\bqab) - a^2(\bqcd) &= (bc+ad+ac)(bc-ad) \label{eq:Id1}\\
    c^2(\bqab) - b^2(\bqcd) &= (ac+bd+bc)(ac-bd) \label{eq:Id2}\\
    d^2(\bqab) - a^2(\bqcd) &= (bd+ac+ad)(bd-ac) \label{eq:Id3}\\
    d^2(\bqab) - b^2(\bqcd) &= (ad+bc+bd)(ad-bc) \label{eq:Id4}
  \end{align}
\end{subequations}

\begin{theorem}
  \label{th:BNumProduct}
  Let $m = \bqab$, $n = \bqcd$ and $k = mn$, with $a, b
  \in \Nonnegative$. 
  \begin{enumerate}
  \item $k = \bqalbet$, with $\alpha, \beta \in \Nonnegative$, has the following solutions:
    \begin{subequations}
      \begin{equation}
        \label{eq:Ident1}
        \begin{cases}
          \alpha = (ad + bc + bd), \quad \beta = (ac - bd), &
          \text{if $ac > bd$;}\\
          \alpha = (ac + ad + bc), \quad \beta = (bd - ac), &
          \text{otherwise.} 
        \end{cases}
      \end{equation}
      \begin{equation}
        \label{eq:Ident2}
        \begin{cases}
          \alpha = (ac + bd + bc), \quad \beta = (ad - bc), &
          \text{if } ad > bc; \\
          \alpha = (ad + ac + bd), \quad \beta = (bc - ad), &
          \text{otherwise.}
        \end{cases}
      \end{equation}
    \end{subequations}

  \item $k = \bqalbet$, with $\alpha, \beta \in \Nonnegative$, has the following solutions:
    \begin{subequations}
      \begin{equation}
        \label{eq:Ident3}     
        \begin{cases}
          \alpha = (ad + bc + bd), \quad \beta = (bd - ac), &
          \text{if $ac < bd$;}\\
          \alpha = (ac + ad + bc), \quad \beta = (ac - bd), &
          \text{otherwise.} 
        \end{cases}
      \end{equation}
      \begin{equation}
        \label{eq:Ident4}
        \begin{cases}
          \alpha = (ac + bd + bc), \quad \beta = (bc - ad), &
          \text{if $ad < bc$;} \\
          \alpha = (ad + ac + bd), \quad \beta = (ad - bc), &
          \text{otherwise.}
        \end{cases}
      \end{equation}    
      \begin{equation}
        \label{eq:Ident5}
        \alpha = (ad + bc + bd), \quad \beta = (ad + bc + ac).
      \end{equation}
      \begin{equation}
        \label{eq:Ident6}
        \alpha = (ac + bd + bc), \quad \beta = (ac + bd + ad).
      \end{equation}
    \end{subequations}
  \end{enumerate}
\end{theorem}

Note that identities~(\ref{eq:Ident3}--\ref{eq:Ident4})
are equivalent to having solutions to $k = \bqalbet$
with non-negative $\alpha$ and non-positive $\beta$ (by
using $-\beta$ instead of $\beta$), and
identities~(\ref{eq:Ident5}--\ref{eq:Ident6}) is
equivalent to having solutions to $k = \bqalbet$ with
non-positive $\alpha$ and $\beta$ (by using $-\alpha$
and $-\beta$ instead of $\alpha$ and $\beta$).

\subsection{Known theorems}
The following well-known theorems are used to prove the
results in this paper.

\begin{theorem}[Fermat]
  \label{th:FLT}
  If $p$ is a prime, and $p \nmid a$, then
  $\congr{a^{p-1}}{0}{p}$.
\end{theorem}
This is known as {\em Fermat's little theorem.} See
~\cite[\S{}6.1]{HAR1995}.
\begin{theorem}[Lagrange]
  \label{th:Lagrange}
Let $p$ be a prime and let $f(x) = a_nx^n + a_{n-1}x^{n-1} +
\cdots + a_0$, where $p \nmid a^n$.  Then the congruence
\congr{f(x)}{0}{p} has at most $n$ distinct solutions
$\alpha$ such that $-p/2 < \alpha \le p/2$.
\end{theorem}
For a proof, see \cite{ALL1995}, Theorem~4.4.1.
\begin{theorem}[Legendre]
  \label{th:KDistinctSolutions}
  If $p$ is a prime and if $k$ divides $(p-1)$, then the
  congruence \congr{x^k-1}{0}{p} 
  has exactly $k$ distinct
  solutions between $-p/2$ and $p/2$.
\end{theorem}
For a proof, see \cite{ALL1995}, Corollary 4.4.6.

\section{Fundamental form of \bprime{}s}
Helguero~\cite{DIC1923},  Fontene~\cite{DIC1923} and Ramanujan~\cite[pp.
259--260]{RAM1988} observed that every prime other than
$3$ that can be represented as \bqab{} is in the form
$6k+1$.  In this section, we prove that it is both the necessary
and sufficient condition, and its
\brep{} is unique (Theorem~\ref{th:BasicFull}).

\begin{theorem}
  \label{th:BasicFull}
  A prime other than $3$ can be represented as
  \bqab{} if and only if it is in the form $6k+1$, and the
  representation is unique.
\end{theorem}

Theorem~\ref{th:BasicFull} states
that 
\begin{enumerate}
\item[(i)] a \bprime{} has a unique
\brep{}, 
\item[(ii)] all \bprime{}s other than $3$ are in the form
$6k+1$, and 
\item[(iii)] all primes in the form $6k+1$ are
\bprime{}s.
\end{enumerate}

(i) is proved in Theorem~\ref{th:PrimeUniqueRep} by assuming that
a \bprime{} has two \brep{}s and deducing that they are
not distinct.  
\begin{theorem}
  A \bprime{} has a unique \brep{}.
  \label{th:PrimeUniqueRep}
\end{theorem}
\begin{proof}
  Let
  \begin{equation}
    \label{eq:PrimeWithTwoRep}
    p = \bqab = \bqcd,
  \end{equation}
  with $a,b,c,d \in \Positive$, be two distinct \brep{}s of
  the prime $p$.
  WLOG, assume $a>c$, so that $ac>bd$.
  Now, using identities \eqref{eq:Id4} and \eqref{eq:Id2},
  \begin{subequations}
    \begin{align}
      p(d^2-b^2) &= (ad-bc)(ad+bc+bd),
      \label{eq:PrimeWithTwoRep2} \\ 
      \intertext{and}
      p(c^2-b^2) &= (ac-bd)(ac+bd+bc). \label{eq:PrimeWithTwoRep3}
    \end{align}
  \end{subequations}

  From \eqref{eq:PrimeWithTwoRep2}, since $p$ is a prime, it should divide at least one of
  $(ad+bc+bd)$ and $(ad-bc)$.  Let us consider each case.

  \begin{description}
  \item [If \boldmath$p|(ad+bc+bd)$\unboldmath] Since $ac>bd$, using \eqref{eq:Ident1},
    $p^2 = \bq{(ad+bc+bd)}{(ac-bd)}$.  Since $ad+bc+bd > 0$, we
    get $p = ad+bc+bd$, giving $ac = bd$.  Now, from
    \eqref{eq:PrimeWithTwoRep3}, $c^2-b^2 = 0$, meaning $c=b$,
    and so, by Theorem~\ref{th:UniqueB}, $a=d$, showing a unique \brep{} of $p$.
    
  \item [If \boldmath$p|(ad-bc)$\unboldmath] Here, we have two cases to consider:
    \begin{description}
    \item [If \boldmath$ad > bc$\unboldmath] By \eqref{eq:Ident2}, $p^2 =
      \bq{(ac+bd+bc)}{(ad-bc)}$, and $p$ should divide
      $(ac+bd+bc)$ as well because $p$ is a prime.  Since
      $(ac+bd+bc) > 0$, this implies $p=(ac+bd+bc)$, meaning
      $ad=bc$.  
    \item [If \boldmath$ad \le bc$\unboldmath] In a similar way, using \eqref{eq:Ident2}, we can show $p =
      (ad+ac+bd)$ and hence $ad = bc$.  
    \end{description}
    We showed that $ad=bc$ in both cases.  Now, from
    \eqref{eq:PrimeWithTwoRep2}, $d^2-b^2 = 0$, meaning $d=b$,
    and so, by Theorem~\ref{th:UniqueB}, $a=c$, showing a unique \brep{} of $p$.
  \end{description}
  We proved that p has a unique \brep{} in both cases.
\end{proof}

(ii) is easily deducible from the properties of
congruences (Theorem~\ref{th:Basic}). 

\begin{theorem}
  \label{th:Basic}
  All \bprime{}s other than $3$ are of the form $6k+1$.
\end{theorem}
\begin{proof}
  Let $p = a^2 + ab + b^2$ be a prime.  
  Let \congr{a}{m}{6}, \congr{b}{n}{6} and
  \congr{\bqab}{z}{6}. Now, 
  \congr{\bq{m}{n}}{z}{6}, using the basic properties
  of congruences.
  For $m = 0 \ldots 5$ and $n = 0 \ldots 5$, $z$ can take
  only the values $0, 1, 3, 4$, i.e., $p$ can take values
  $6k$, $6k+1$, $6k+3$ and $6k+4$.  Here $6k$ and $6k+4$ are
  always composite.  $6k+3$ is composite except for $k
  = 0$, i.e., when $p = 3$.  So, the only prime values $p$ can
  take are $3$ and $6k+1$.
\end{proof}

(iii) is proved in four steps:
\begin{enumerate}
\item We show that if a \bnum{} is divided by a \bprime{},
  we get another \bnum{}. (Theorem~\ref{th:BNumByBPrime})
\item Using the previous result, we show that if we divide a
  \bnum{} with a factor that is not a \bnum{}, if such a factor
  exists, then at least one prime factor of the quotient is
  not a \bprime. (Theorem~\ref{th:BNumFactorNonBNum})
\item Now we show that every factor of a \bnum{} \bqab{},
  with $a$ and $b$ being relatively prime, is a \bnum{}.  We start
  with the assumption that there is a factor that is not a
  \bnum{}.  We then prove that for any number with that
  property, we can find a smaller positive number with the
  same property.  Now, by the principle of infinite descent,
  there is no such number. (Theorem~\ref{th:FactorOfBNum})
\item Finally, using some well-known results, we prove that
  every prime in the form $6k+1$ divides a number $a^2+a+1$
  for some $a \in \Positive$.  Since $a^2+a+1$ is a \bnum{} and $a$ and
  $1$ are coprimes, the previous result implies that the
  prime is a \bprime. (Theorem~\ref{th:BasicConv})
\end{enumerate}

\begin{theorem}
  \label{th:BNumByBPrime}
  If a \bnum{} $n$ has a \bprime{} factor $p$, then
  $(n/p)$ is a \bnum.
\end{theorem}

\begin{proof}
  Let
  \begin{equation}
    \label{eq:NEqualsBnum}
    n=\bqab,
  \end{equation}
  and $n = p \cdot q$, where $p =
  \bqcd$, with $a, b, c, d \in \Nonnegative$.  We need to prove that $q$ is a \bnum.

  Consider the identity~\eqref{eq:Id1}.
  Since $p$ divides the LHS, it should divide the RHS.  Since $p$ is a
  prime, it divides at
  least one of $(bc+ad+ac)$ and $(bc-ad)$ .  That leads
  to the following two cases.
  \begin{description}
  \item[Case 1] $p$ divides $(bc+ad+ac)$.

    Let $(bc+ad+ac) = rp$, with $r \in \Integer$.
    Set $a=rd+y$ and $b=rc+x$, where $x,y \in \Integer$. 
    Combining these, $r (\bqcd) + cx + dy + cy =
    rp$, which means
    \begin{equation}
      cx + cy + dy = 0.
      \label{eq:CxCyDySumZero}
    \end{equation}
    Now, substituting $a=rd+y$ and $b=rc+x$ in
    \eqref{eq:NEqualsBnum} and using \eqref{eq:CxCyDySumZero}, we get
    \begin{equation}
      n =  r^2(\bqcd) + \bqxy + r(cx+dx+dy).
      \label{eq:NsRep1}
    \end{equation}
    Since \eqref{eq:CxCyDySumZero} can be rewritten as
    \begin{equation*}
      c(x+y) + dy = 0,
    \end{equation*}
    and $(c,d) = 1$, $c$ divides $y$.  Let $y =
    -cw$, so that $x = (c+d)w$, where $w \in \Integer$.
    Substituting these values of $x$ and $y$ in
    \eqref{eq:NsRep1} and simplifying, we get
    \begin{equation*}
      n = (\bqcd)(\bq{r}{w}),
    \end{equation*}
    which means $q = \bq{r}{w}$.
  \item[Case 2] $p$ divides
    $(bc-ad)$.

    Let $(bc-ad) = rp$, where $r \in \Integer$. Set
    $a=-rd+y$ and $b=rc+x$,  where $x,y \in \Integer$.
    Combining these,   $rc^2 + cx - dy + rd^2 = rc^2
    + rcd + rd^2$, which means
    \begin{equation}
      \label{eq:RcdPYdMXcE0}
      rcd + dy - cx = 0.
    \end{equation}
    Now, substituting $a=y-rd$ and $b=x+rc$ in
    \eqref{eq:NEqualsBnum} and using
    \eqref{eq:RcdPYdMXcE0}, we get
    \begin{equation}
      \label{eq:NsRep2}
      \begin{split}
        n = r^2(\bqcd) + (\bqxy) +  r(cy - dx).
      \end{split}
    \end{equation}
    Since \eqref{eq:RcdPYdMXcE0} can be rewritten as
    \begin{equation*}
      c(x - rd) = dy,
    \end{equation*}
    and $(c,d) = 1$, $c$ divides $y$. Let $y = cw$, so that $x = (r+w)d$,
    where $w \in \Integer$.
    Substituting these values of $x$ and $y$ in
    \eqref{eq:NsRep2} and simplifying, we get
    \begin{equation*}
        n = (\bqcd)(\bq{r}{w}),
    \end{equation*}
    which means $q = \bq{r}{w}$.
  \end{description}
  So, in either case, $q = \bq{r}{w}$, where $r, q \in
  \Integer$.  So, by Theorem~\ref{th:BNumOfAllIntegers}, it is a \bnum.
\end{proof}

\begin{theorem}
  \label{th:BNumFactorNonBNum}
  If a \bnum{} $n$ has a factor $m$ which is not a \bnum, then
  $(n/m)$ has at least one prime factor that is not a \bprime.
\end{theorem}

\begin{proof}
  Let $n = m \cdot k$.  Factor $k$ into prime factors $k =
  p_1 \cdot p_2 \cdots p_x$.  Suppose that all $p_i$s ($i = 1, 2, \cdots
  x$) are \bprime{}s.  Now, by Theorem~\ref{th:BNumByBPrime},
  $n/p_1$ is a \bnum, and hence, $n/(p_1 \cdot p_2)$ is a \bnum,
  and continuing this, $n/(p_1 \cdot p_2 \cdots p_x)$ is a \bnum,
  which means $m$ is a \bnum, which is a contradiction, so all
  $p_i$s cannot be \bprime{}s.
\end{proof}

\begin{theorem}
  \label{th:FactorOfBNum}
  If $n = \bqab$ where $(a, b) = 1$, then each factor of
  $n$ is a \bnum.
\end{theorem}

\begin{proof}
  Suppose that $n$ has a factor $p$ that is not a \bnum.  Let $a =
  xp + \alpha$ and $b = yp + \beta$, with $\alpha,\beta\in\Integer$,
  so that $-p/2 < \alpha, \beta \le p/2$.  Since $p$ divides
  $(\bqab)$, it should divide
  $(\bq{\alpha}{\beta})$ as well.  Let $\bq{\alpha}{\beta} = pq$.
  Since $(\bq{\alpha}{\beta}) \ge 0$, because of Theorem~\ref{th:BqNotNegative},
  $pq \ge 0$, and hence $q \ge 0$.  Also, since
  $(\bq{\alpha}{\beta}) \le 3p^2/4$, $q \le 3p/4 < p$.

  Let $(\alpha, \beta) = \gamma$.  Also let $u =
  (\alpha/\gamma)$ and $v = (\beta/\gamma)$.  Now 
  \begin{equation*}
    \bq{u}{v} = \frac{pq}{\gamma^2}. 
  \end{equation*}
  Since $\gamma$ does not divide $p$ (otherwise $\gamma$
  will divide both $a$
  and $b$), $\gamma^2$ should divide $q$.  Let $k=
  (q/\gamma^2)$.  Now,
  \begin{equation*}
    \bq{u}{v} = pk.
  \end{equation*}
  Since $p$ is not a \bnum, because of
  Theorem~\ref{th:BNumFactorNonBNum}, $k$ must have a prime
  factor, say $r$, 
  that is not a \bprime.  Since $r \le k \le q < p$, $r <
  p$. 

  So, we started with a number $p \ge 0$, which is not a
  \bnum{} but is a factor of a \bnum{}, and found a
  smaller number $r \ge 0$ with the same property.  So, by the
  principle of infinite descent, this is impossible.
  Hence $p$ must be a \bnum.
\end{proof}

If $(a,b) = k$, $\bqab = k^2x^2 + kx \cdot ky + k^2y^2 = k^2
\cdot (\bqab)$.  So, if a \bnum{} \bqab{} is square-free,
$(a,b) = 1$.  This leads to the following corollary.

\begin{corollary}
  A square-free \bnum{} is the product of \bprime{}s.
  \label{th:SquareFreeBNumIsBPrimeProduct}
\end{corollary}

\begin{theorem}
  For every prime $p$ of the form $6k+1$, there exists a
  unique positive integer $z < p/2$, such that
  \congr{\bqz{z}}{0}{p}.
  \label{th:bqz}
\end{theorem}
\begin{proof}
  By Fermat's Little Theorem (Theorem~\ref{th:FLT}), 
  \begin{equation}
    \label{eq:BasicConv1}
    \congr{x^{p-1} - 1}{0}{p}.
  \end{equation}
  Since $(p-1)/2 = 3k$, $x^{p-1} - 1 = (x^{3k} + 1)(x^{3k}-1)$,
  and the solutions of \eqref{eq:BasicConv1} are given by the
  solutions of 
  \begin{subequations}
    \begin{align}
      x^{3k}+1 &\equiv 0 \pmod{p},\label{eq:BasicConv2}
      \intertext{and} 
      x^{3k}-1 &\equiv 0 \pmod{p}.\label{eq:BasicConv3}
    \end{align}
  \end{subequations}
  By Theorem~\ref{th:KDistinctSolutions}, \eqref{eq:BasicConv3} has
  exactly $3k$ solutions between $-p/2$ and $p/2$.  Now,
  the solutions of \eqref{eq:BasicConv3} are given by the
  solutions of 
  \begin{subequations}
    \begin{align}
      x^k-1 &\equiv 0 \pmod{p},\label{eq:BasicConv4}
      \intertext{and}
      x^{2k}+x^k+1 &\equiv 0 \pmod{p}. \label{eq:BasicConv5}
    \end{align}
  \end{subequations}

  \eqref{eq:BasicConv4} has exactly $k$ solutions such
  that $-p/2 < x < p/2$, so \eqref{eq:BasicConv5}
  should have $2k$ solutions.  For any of these $2k$
  solutions, $y=x^k$ gives a solution to 
  \begin{equation}
    \label{eq:BasicConv6}
    \congr{\bqz{y}}{0}{p}.
  \end{equation}
  Now, any $w$ such that
  $w=y+k\cdot p$, where $k \in \Integer$, also satisfy \eqref{eq:BasicConv6}.
  So, we can find a $u$ such that
  $-p/2 < u<p/2$ that satisfies
  $\congr{\bqz{u}}{0}{p}$.  Now, setting $v=-(u+1)$,
  we find $\bqz{v}=\bqz{u}$, so $\congr{\bqz{v}}{0}{p}$.  If
  $u<0, 0\le v < p/2$, and vice versa; so there are at
  least two
  solutions, greater than  $-p/2$ and less
  than $p/2$, of which one is negative, and the other
  is not.    But by Lagrange's
  theorem (Theorem~\ref{th:Lagrange}), \eqref{eq:BasicConv6} has at most $2$
  solutions for  $-p/2 < y < p/2$.  Since $p$ is odd
  and $p>1$, we can conclude that \eqref{eq:BasicConv6} has exactly
  two solutions between $-p/2$ and $p/2$, one a positive integer less than
  $p/2$, and the other a negative integer greater than
  $-p/2$.
 \end{proof}

\begin{theorem}
  \label{th:BasicConv}
  If a prime number is of the form $6k+1$, it is a \bprime.
\end{theorem}
\begin{proof}
  Let $p = 6k+1$ be a prime.  By Theorem~\ref{th:bqz}, there
  exists a $z \in \Positive$  such that
  \congr{\bqz{z}}{0}{p}.  This means $\bqz{z} = m \cdot
  p$, where $m \in \Positive$.  Since $(z,1) = 1$, by Theorem~\ref{th:FactorOfBNum}, every factor of
  \bqz{z} should be a \bnum{}, so $p$ is a \bprime.
\end{proof}

Proof of Theorem~\ref{th:BasicFull} follows from
Theorems \ref{th:PrimeUniqueRep}, \ref{th:Basic} and
~\ref{th:BasicConv}.

\section{Factors of \bnum{}s}
In this section, we prove the following theorem.

\begin{theorem}
  \label{th:BNumFactorization}
  The necessary and sufficient condition of any
  non-negative integer to be in the form \bqab{} is
  that, in its prime factorization, all primes other
  than $3$ that are
  not in the form $(6k+1)$ have even exponents.
\end{theorem}

The {\em sufficient} part is proved first:

\begin{theorem}
  \label{th:BNumFactorization1}
  If the factorization of a number has even exponents
  for all primes other than $3$ and those in the form
  $(6k+1)$, the number is a \bnum.
\end{theorem}
\begin{proof}
  From \eqref{eq:Ident1}, the product of any numbers
    of \bnum{}s is a \bnum{}.  In particular, the product of
    any number of \bprime{}s is a \bnum.
  Since every square is a \bnum{}, multiplying a
    \bnum{} with a square will yield a \bnum.
  Hence, the product of any number of \bprime{}s and
    any square is a \bnum.

The theorem easily follows from this.
\end{proof}

We prove the {\em necessary} part as follows:

\begin{enumerate}
\item We prove that if we divide a \bnum{} by a square
  factor, if it has one, the quotient is a \bnum.(Theorem~\ref{th:BNumDividedBySquare})
\item Now, using
  Corollary~\ref{th:SquareFreeBNumIsBPrimeProduct}, the
  result follows.
\end{enumerate}

\begin{theorem}
  If the product of $m$ and $k^2$, with $m,k \in \Positive$, is a \bnum, then $m$ is a \bnum.
  \label{th:BNumDividedBySquare}
\end{theorem}
\begin{proof}
  We have 
  \begin{equation}
    \label{eq:FcnContainsSquare}
    \bq{a}{b} = k^2 \cdot m.
  \end{equation}
  Let $s^2$ be the largest square factor of $m$.  Let
  $n=m/s^2$. We need to only show that $n$ is a \bnum{},
  because, $m$, being the product of two \bnum{}s, is a \bnum{}
  then.  Now, defining $p=ks$, \eqref{eq:FcnContainsSquare}
  becomes
  \begin{equation*}
    \bq{a}{b} = p^2 \cdot n,
  \end{equation*}
  where $n$ does not have square factor $>1$.  Let $g=(a,b)$
  and $c=a/g, d=b/g$. Now $g^2(\bqcd) = p^2n$, or
  \begin{equation*}
    \bqcd = \frac{p^2}{g^2} \cdot n.
  \end{equation*}
  Since $g^2$ doesn't divide $n$ if $g>1$, $p/g$ must be an integer.
  Now, by Theorem~\ref{th:FactorOfBNum}, each factor of the LHS must
  be a \bnum{}, so each factor of $n$ must be a \bnum.  So,
  $n$, being the product of \bnum{}s, must be a \bnum{}. 
\end{proof}

\begin{theorem}
  \label{th:BNumFactorization2}
  Every prime factor, which is not a \bprime,  of a \bnum{}  has 
  even exponents in the standard form of expansion.
\end{theorem}
\begin{proof}
  Let $n$ be any \bnum{} and let $s^2$ be its largest square
  factor.  By Theorem~\ref{th:BNumDividedBySquare},
  $(n/s^2)$ is a square-free \bnum{}, which is, by
  Corollary~\ref{th:SquareFreeBNumIsBPrimeProduct}, the
  product of \bprime{}s.  So, all non-\bprime{} factors
  should be contained in $s^2$, hence 
  should have even exponents.
\end{proof}

Proof of Theorem~\ref{th:BNumFactorization} follows from Theorems
\ref{th:BNumFactorization1} and \ref{th:BNumFactorization2}.

\section{General form of a \bnum}
From the already proved theorems, the general form of a
\bnum{} is
\begin{equation}
  \label{eq:BNumGenForm}
  n = x^2 \cdot 3^y \cdot a^\alpha{}\cdot{}b^\beta{}\cdot{}c^\gamma\ldots
\end{equation}
where

\begin{tabular}{lcp{4in}}
  $x$ &:& Some number in \Nonnegative, with no
  prime factor in the form $6k+1, k \in \Positive$. \\
  $y$ &:& Some number in \Nonnegative \\
  $a, b, c,\ldots$ &:& A prime in the form $6k+1, k \in \Nonnegative$ \\
  $\alpha, \beta, \gamma, \ldots$ &:& Some number in \Nonnegative 
\end{tabular}

\vspace{10pt}

A number is a \bnum{} if and only if it can be represented
in the form \eqref{eq:BNumGenForm}.
\begin{conjecture}
  \label{th:BNumNumberOfRep}
  The number of distinct
  \brep{}s of a \bnum{} represented in the form given by \eqref{eq:BNumGenForm}, counting the cases when $a=b$ and when
  either of $a$ or $b$ is zero, is
  given by
  \begin{equation}
    \begin{cases}
      \frac{1}{2} \Big(1 + (\alpha+1)(\beta+1)(\gamma+1)\ldots\Big), & \quad
      \text{if all of } \alpha, \beta, \gamma, \ldots \text{ are even;} \\
      \frac{1}{2} (\alpha+1)(\beta+1)(\gamma+1)\ldots, & \quad
      \text{otherwise.}
    \end{cases}
    \label{eq:NumOfRep}
  \end{equation}
\end{conjecture}

The reader will recognize that \eqref{eq:NumOfRep} is the
same as the expression for the number of distinct
representations of a number as a sum of two squares, a
result given by Gauss and Legendre, later proved by Jacobi
and others~\cite{DIC1920}.  This expression applies here also,
with a different definition of the individual parameters.  A
proof is not attempted in this paper.

\section{\bnum{}s in terms of rational numbers}
Theorem~\ref{th:RationalBNum} 
extends the domain of $a$ and $b$ in the previous theorems
from \Nonnegative{} to \NonnegativeRational.

\begin{theorem}
\label{th:RationalBNum}
If a number $n$ is representable as \bq{\alpha}{\beta},
where $\alpha$ and $\beta$ are positive rational numbers,
then $n$ is a \bnum. 
\end{theorem}
\begin{proof}
Let $\alpha = (a/b)$ and $\beta = (c/d)$.  Now
  \begin{equation*}
    \begin{split}
      n &= \bq{\alpha}{\beta} \\
      &= \bq{\left(\frac{a}{b}\right)}{\left(\frac{c}{d}\right)} \\
      &= \frac{a^2d^2 + abcd + b^2c^2}{b^2d^2}
    \end{split}  
  \end{equation*}
Hence,
\begin{equation*}
  n \cdot (bd)^2 = \bq{(ad)}{(bc)}
\end{equation*}
Then, by Theorem~\ref{th:BNumDividedBySquare}, $n$ is a \bnum.
\end{proof}

\section{Comparison with the form \bqmab}

Being another binary quadratic form with the same
discriminant, the form \bqmab{} shares many
properties with \bqab. 

All primes in the form $6k+1$ and $3$ can be represented as
\bqmab, but this representation is not unique.  In general,
if $n = \bqab$, where $n \in \Nonnegative$, then $n$ can be represented as
\bqmxy{} in two ways by setting $ x = a,  y = a + b$
and $x = b, y = a + b$.  The unique representation occurs only when $a=b$, i.e., when
$n=3a^2$, and $3$ is the only prime that has unique
representation.  Similarly, when $n=\bqmab$, with $a \ge b$, $n$ can be
represented as \bqxy{} as well, by setting $x = a-b, y =
b$. Because of this equivalence, all results in this paper,
except Theorem~\ref{th:BasicFull},
Theorem~\ref{th:PrimeUniqueRep} 
 and \eqref{eq:NumOfRep}, are
valid for the form \bqmab{} as well.

\section{Comparison with forms equivalent to $ka^2+kab+kb^2$}
Many results discussed in this paper have been proved for
some of the binary quadratic forms equivalent to $ka^2+kab+kb^2, k \in
\Positive$.  For example, Euler (1763) proved the following for the form $a^2+3b^2$
~\cite{DIC1923}, which is equivalent to $2a^2+2ab+2b^2$.
\begin{enumerate}
\item A prime $p$ can be represented in this form if and only if
  $p=3$ or \congr{p}{1}{6} (cf. Theorem~\ref{th:BasicFull}).
  Goring (1874) proved the
  uniqueness of this representation (cf. Theorem~\ref{th:PrimeUniqueRep}).
\item Every prime divisor (other than $2$) of a number in
  this form, with $a$ and $b$ coprimes, is itself in this
  form (cf. Theorem~\ref{th:FactorOfBNum}).
\end{enumerate}

Similar analogies can be found for forms equivalent to any
$ka^2+kab+kb^2$, like $a^2+ab+7b^2$, $a^2+12b^2$,
$3a^2+4b^2$, $a^2+ab+19y^2$,
$3a^2+3ab+7b^2$, etc.

\section{Comparison with the form $a^2+b^2$}
Every theorem in this paper has an analogous theorem
for the form corresponding to the sum of two squares.
They are listed below:
\begin{enumerate}
\item (L. Pisano, 1225) The product of two numbers can be expressed as
  the sum of two squares if the individual numbers can
  be.  \cite[Ch. 6]{DIC1920} (cf. Theorem~\ref{th:BNumProduct}.)
\item (Euler, 1749) A prime $p$ can be expressed as the sum of two
  squares if and only if $\congr{p}{1}{4}$, with the
  only exception of $2$. (cf. Theorem~\ref{th:BasicFull}.)  Gauss
  proved that this representation is unique. (cf. Theorem~\ref{th:PrimeUniqueRep})
\item (Euler, 1749) If $m=a^2+b^2=pl$, where
  $p=c^2+d^2$ is a prime, then $l$ is the sum of two
  squares. \cite[\S{}7.1]{ALL1995} (cf. Theorem~\ref{th:BNumByBPrime})
\item (Euler, 1749) Let $m=a^2+b^2=nk$. If $k$ is not a
 sum of two squares, then $n$ has a prime factor which
 is not a sum of two squares. \cite[\S{}7.1]{ALL1995}
 (cf. Theorem~\ref{th:BNumFactorNonBNum}.)
\item (Euler) A positive number is the sum of two squares if
  and only if all of its factors in the form $4m+3$
  have even exponents in the standard
  factorization. \cite[Th. 366]{HAR1995} (cf. Theorem~\ref{th:BNumFactorization}.)
\item (Goldbach, 1743) A prime $4k-1$ cannot divide the sum of two
  relatively prime squares. \cite[Ch. 6]{DIC1920}
  (cf. Theorem~\ref{th:FactorOfBNum}.)
\item $-1$ is a quadratic residue of primes of the form
  $4k+1$ \cite[Th. 82]{HAR1995}.  This means there exists a
  positive integer $z$ such that \congr{z^2+1}{0}{p} if
  $p=4k+1$ is a prime. (cf. Theorem~\ref{th:bqz}.)
\item If $n = \alpha^2 + \beta^2$, where $\alpha$ and
  $\beta$ are rational, then $n$ is the sum of two
  integer squares.  \cite[Note 7.1.5]{ALL1995} (cf. Theorem~\ref{th:RationalBNum}.)
\item (Gauss, 1801) The general form of a number $n$ that can be
  expressed as the sum of two squares is
  \begin{equation}
  n = x^22^ya^\alpha{}b^\beta{}c^\gamma\ldots
  \label{eq:NumRepSumOfSq}
  \end{equation}
  where
  $a, b, c, \ldots$ are primes of the form
  $4k+1$, and $x$ is the product of primes of the form
  $4k+3$. (cf. \eqref{eq:BNumGenForm}.)
\item Jacobi (1834) gave an expression for the number of
  representations of a number $n$ in the form $a^2+b^2$, which,
  when applied to distinct representations with $a$ and $b$ in \Nonnegative,
  reduces to \eqref{eq:NumOfRep}, when $n$ is expressed
  as in \eqref{eq:NumRepSumOfSq}.  (cf. \eqref{eq:BNumGenForm} and
  \eqref{eq:NumOfRep}.)
\end{enumerate}

\section{Conclusions}
Many mathematicians have observed the properties of the form
\bqab{} and similar forms, and quite an extensive study of
binary quadratic forms has been done in the past.  Many of
the results in this paper can be obtained using alternate
techniques.  One such method is to assume the already known
properties of the form $a^2+3b^2$, and then deduce the
result for \bqab{} by trivial substitutions.

The contribution of this paper is the proof to these
properties by elementary theory of numbers.  It also
conjectures a formula for the number of representations of a
number in the form.

This paper proves that the binary quadratic form \bqab{} has
properties very much analogous to those of the form
$a^2+b^2$.   It is possible that these are two special cases
of a more general form.  More investigation is needed in
this area.

\providecommand{\bysame}{\leavevmode\hbox to3em{\hrulefill}\thinspace}
\providecommand{\MR}{\relax\ifhmode\unskip\space\fi MR }
\providecommand{\MRhref}[2]{%
  \href{http://www.ams.org/mathscinet-getitem?mr=#1}{#2}
}
\providecommand{\href}[2]{#2}

\end{document}